\documentclass[12pt]{amsart}
\usepackage{amsmath,amssymb,amsbsy,amsfonts,latexsym,amsopn,amstext,cite,
                                               amsxtra,euscript,amscd,bm}
\usepackage{url}

\usepackage{mathrsfs}

\usepackage{color}
\usepackage[colorlinks,linkcolor=blue,anchorcolor=blue,citecolor=blue,backref=page]{hyperref}
\usepackage{color}
\usepackage{graphics,epsfig}
\usepackage{graphicx}
\usepackage{float}
\usepackage{epstopdf}
\hypersetup{breaklinks=true}

\usepackage[np]{numprint}
\npdecimalsign{\ensuremath{.}}

\usepackage{bibentry}
\usepackage[capitalize]{cleveref}

\usepackage[english]{babel}
\usepackage{mathtools}
\usepackage{todonotes}
\usepackage{url}
\usepackage{enumitem}

\newtheorem{theorem}{Theorem}

\DeclareMathOperator{\Sym}{Sym}
\DeclareMathOperator{\GL}{GL}
\DeclareMathOperator{\SL}{SL}
\DeclareMathOperator{\SU}{SU}
\DeclareMathOperator{\GU}{GU}
\DeclareMathOperator{\Sp}{Sp}

\usepackage{xcolor}

\newif\ifcomment

\newcommand{\F}{\mathbb{F}}
\newcommand{\inmod}{\hspace{-.8ex}\mod}
\newcommand{\dispmod}{\hspace{-1.7ex}\mod}

\let\cal\mathcal
\let\tensor\otimes
\let\eps\varepsilon
\let\leq\leqslant
\let\geq\geqslant

\allowdisplaybreaks

\begin{document}

\title[A faster method to construct extraspecial normaliser subgroups]{A faster method to construct extraspecial normaliser subgroups}

\author[M. Hauck] {Maximilian Hauck}
\address{Max Planck Institute for Mathematics, Vivatsgasse 7, 53111 Bonn, Germany}
\email{max.hauck01@gmail.com}

\begin{abstract}    
We show how to improve the runtime of the construction of generators of maximal subgroups of $\SL(d, q), \SU(d, q)$ and $\Sp(d, q)$ which arise as normalisers of extraspecial groups or $2$-groups of symplectic type given in \cite[§9]{HoltRoneyDougal} from $O(d^3\log d\log q+\log^2 q)$ to $O(d^2\log d\log^{1+\eps} q+\log^{2+\eps} q)$.
\end{abstract}

\keywords{Finite classical groups, maximal subgroups, Aschbacher classes, polynomial time algorithms}
\subjclass[2020]{20D06, 20E28, 20G40, 68Q25}

\maketitle

\section{Introduction}

\subsection{Motivation and background}  
In \cite{CannonHolt} and \cite{EickHulpke}, it has been shown that the computation of the maximal subgroups of an arbitrary finite permutation group $G$ effectively reduces to the case where $G$ is almost simple; here, the case where $G$ is a finite classical group is of particular importance. Consequently, it is a problem of major interest in computational group theory to construct the maximal subgroups of classical groups.

These maximal subgroups of a finite classical group $G$ can be partitioned into nine classes according to Aschbacher's theorem, see \cite{Aschbacher}: The first eight of these classes are usually denoted $\cal{C}_1$ through to $\cal{C}_8$; they as well as the structure and conjugacy of their members and their maximal subgroups are described in \cite{KleidmanLiebeck}. The last class, usually denoted by $\cal{S}$, consists of almost simple groups acting absolutely irreducibly and satisfying some additional conditions, see \cite[§1.2]{KleidmanLiebeck}.

For the classes $\cal{C}_1$ through to $\cal{C}_8$, polynomial time algorithms to construct all these maximal subgroups in terms of generators (up to conjugation by the elements of a certain subgroup of the general linear group) have been given by Holt and Roney-Dougal: In \cite{HoltRoneyDougal}, they treat the case of linear, symplectic and unitary groups while \cite{HoltRoneyDougal2} is devoted to the case of orthogonal groups. In this note, we want to give an improvement upon their construction of the maximal subgroups in class $\mathcal{C}_6$ for $G$ a linear, symplectic or unitary group which reduces the time complexity of the necessary computations from $O(d^3\log d\log q+\log^2 q)$ to $O(d^2\log d\log^{1+\eps} q+\log^{2+\eps} q)$; here, $d$ is the dimension of the tautological representation of $G$ and $q$ is the cardinality of the finite field we are working over.

\subsection{Results}

Let $q=p^e$ be a prime power, $r$ a prime dividing $q-1$ and $d=r^m$ for some integer $m\geq 1$. Moreover, let $G=\SL(d, q), G=\SU(d, q)$ or $G=\Sp(d, q)$ and we denote these three cases by $\mathsf{L}, \mathsf{U}$ and $\mathsf{S}$, respectively.

There are two isomorphism types of extraspecial groups of order $r^{1+2m}$, see \cite[Thm. 5.2]{Gorenstein}. For odd $r$, we only consider the isomorphism type which has exponent $r$ and is denoted $r_+^{1+2m}$ since the normaliser in $\GL(d, q)$ of an extraspecial group of the other isomorphism type is properly contained in the normaliser of an extraspecial group of type $r_+^{1+2m}$. For $r=2$, the extraspecial groups of minus type are central products of a quaternion group of order $8$ with any number of dihedral groups of order $8$; taking a central product of such an extraspecial group of isomorphism type $2_-^{1+2m}$ with a cyclic group of order $4$ yields a $2$-group of symplectic type. Under the condition $r\mid q-1$ we have imposed, both the extraspecial groups of order $r^{1+2m}$ and the symplectic type groups of order $2^{2+2m}$ act on the vector space $\F_q^d$ and this action is absolutely irreducible.

From \cite[Table 4.6.B]{KleidmanLiebeck}, one can extract the conditions under which the normalisers of one of the aforementioned extraspecial or symplectic type groups in $G$ are maximal subgroups of $G$. We collect these results in \cref{table:maximals}; here, the types 
$r_+^{1+2m}.\Sp(2m, r), (4\circ 2_-^{1+2m}).\Sp(2m, 2)$ and $2_-^{1+2m}.O^-(2m, 2)$ indicate whether the corresponding subgroup is the normaliser in $G$ of an extraspecial group of plus type and order $r^{1+2m}$, of a $2$-group of symplectic type of order $2^{2+2m}$ or of an extraspecial group of minus type of order $2^{1+2m}$, respectively, as described in the previous paragraph. In all cases, it is additionally necessary to assume that $e$ is minimal with respect to the condition $r^*\mid q-1=p^e-1$, where $r^*=r$ if $r$ is odd and $r^*=4$ if $r=2$.

\begin{table}[h!]
\centering
\begin{tabular}{c|c|c}
Type & Case & Conditions \\ [.5ex]
\hline && \\ [-1.5ex]
$r_+^{1+2m}.\Sp(2m, r)$ & $\mathsf{L}$ & $er$ odd \\
 & $\mathsf{U}$ & $e$ even, $r$ odd \\
$(4\circ 2_-^{1+2m}).\Sp(2m, 2)$ & $\mathsf{L}$ & $e=1$, $m\geq 2$ \\
 & $\mathsf{U}$ & $e=2$ \\
$2_-^{1+2m}.O^-(2m, 2)$ & $\mathsf{L}$ & $e=1, m=1$ \\
 & $\mathsf{S}$ & $e=1$
\end{tabular}
\vspace{1em}
\caption{Maximal subgroups of $G$ in class $\mathcal{C}_6$}
\label{table:maximals}
\end{table}

Improving upon \cite[Prop. 9.5]{HoltRoneyDougal}, our main result is the following:

\begin{theorem}
\label{thm:main}
In each of the cases $\mathsf{L}, \mathsf{U}$ and $\mathsf{S}$, one can construct generators for representatives of the conjugacy classes in $\GL(d, q)$, $\GU(d, q)$ or $\Sp(d, q)$, respectively, of all groups in \cref{table:maximals} in time $O(d^2\log d\log^{1+\eps} q+\log^{2+\eps} q)$.
\end{theorem}

\subsection{Notations and conventions}

As in the previous section, throughout the paper, let $q=p^e$ be a prime power, $r$ a prime dividing $q-1$ and $d=r^m$ for some integer $m\geq 1$; if $r=2$, we impose the stronger assumption that $4$ divides $q-1$. Moreover, we fix a primitive $r$-th root of unity $\omega$ in $\F_q$. The group $G$ always denotes one of $\SL(d, q), \SU(d, q)$ or $\Sp(d, q)$ and we denote these cases by $\mathsf{L}, \mathsf{U}$ and $\mathsf{S}$, respectively; note that our convention for unitary groups is that $\SU(d, q)\leq \GL(d, q)$ and hence, in particular, the group $\SU(d, q)$ is only defined if $q$ is a square.

We use $I_d$ to denote the identity matrix of size $d\times d$. Furthermore, for integers $a, n\geq 0$, we use $a\inmod n$ to denote the remainder of $a$ upon division by $n$, i.e.\ the unique integer $0\leq b\leq n-1$ such that $a-b$ is divisible by $n$. 

For an $s\times t$-matrix $A$ and a $u\times v$-matrix $B$, we define the Kronecker product $A\tensor B$ as the $su\times tv$-matrix whose $((i-1)u+k, (j-1)v+l)$ entry is $A_{ij}B_{kl}$ for $1\leq i\leq s, 1\leq j\leq t, 1\leq k\leq u$ and $1\leq l\leq v$. The Kronecker product is associative and satisfies $(A\tensor B)(C\tensor D)=AC\tensor BD$ for any matrices $A, B, C$ and $D$ such that the products $AC$ and $BD$ are defined. Most importantly, for square matrices $A$ and $B$ of size $s\times s$ and $u\times u$, respectively, we have $\det(A\tensor B)=\det(A)^u\det(B)^s$.

We assume that primitive polynomials and primitive field elements of $\F_q$ are known and that elements of $\F_q$ are represented by polynomials over $\F_p$ of degree $e-1$. Consequently, by \cite[Thm. 8.23]{GathenGerhard}, field operations in $\F_q$ require time $O(\log^{1+\eps} q)$ and, in particular, computing the element $\omega$ above takes time $O(\log^{2+\eps} q)$ using exponentiation by squaring. We note that, here and throughout the paper, the constant $\eps>0$ is arbitrary and does not affect the implicit constant in the $O$ notation.

\section{Proof of \cref{thm:main}} 

We now show how to modify the construction from \cite[Prop. 9.5]{HoltRoneyDougal} to obtain the improved time complexity estimate of \cref{thm:main}. For $d\leq 4$, the construction from loc.\ cit.\ already achieves the desired time complexity and thus we may without loss of generality assume that $d=r^m>4$. In the following, let $E$ be an extraspecial group or a $2$-group of symplectic type subject to the conditions of \cref{table:maximals}.

Now recall that, in loc.\ cit., the authors construct matrices $X_i, Y_i, U_i, V_i$ and $W_i$ for $1\leq i\leq m$ which, together with scalar matrices, generate the normaliser of $E$ in $\GL(d, q)$; scaling these matrices so that their determinant is $1$ then yields a set of generators of the maximal subgroup of $G$ described by \cref{table:maximals}. By \cite[Lem.s 9.2, 9.3, 9.4]{HoltRoneyDougal}, the matrices $X_i$, $Y_i$, $U_i$, $V_i$ and $W_i$ can be constructed in time $O(d^2\log d\log^{1+\eps} q+\log^{2+\eps} q)$ and thus it suffices to show how to carry out the rescaling step in the desired time.

\subsection{Maximal subgroups of type $r_+^{1+2m}.\Sp(2m, r)$}

We start with the case where $E$ is an extraspecial group of isomorphism type $r_+^{1+2m}$ for $r$ odd. For this, recall the definitions of the following matrices from \cite[Lem.s 9.1, 9.2]{HoltRoneyDougal}: 
\begin{enumerate}[label=--]
\item The matrix $X\in\GL(r, q)$ is diagonal with $X_{ii}=\omega^{i-1}$ for $1\leq i\leq r$.
\item The matrix $Y\in\GL(r, q)$ is the permutation matrix defined by the cyclic permutation $(1, 2, \dots, r)$.
\item The matrix $U\in\GL(r, q)$ is diagonal with $U_{ii}=\omega^{i(i-1)/2}$ for $1\leq i\leq r$.
\item The matrix $V\in\GL(r, q)$ has entries $V_{ij}=\omega^{(i-1)(j-1)}$ for $1\leq i, j\leq r$.
\item The matrix $W\in\GL(r^2, q)$ is the permutation matrix defined by the permutation $w\in\Sym(\{0, \dots, r^2-1\})$ given by $a\mapsto (a+((a-1)\inmod r)\cdot r)\inmod r^2$.
\end{enumerate}
Then the matrices $X_i, Y_i, U_i$ and $V_i$ are defined by $X_i\coloneqq I_{r^{m-i}}\tensor X\tensor I_{r^{i-1}}$ for $1\leq i\leq m$ and similarly for $Y_i, U_i$ and $V_i$. Analogously, for $m\geq 2$, we define $W_i\coloneqq I_{r^{m-1-i}}\tensor W\tensor I_{r^{i-1}}$ for $1\leq i\leq m-1$. Note that this means 
\begin{equation}
\label{eq:detXi}
\det(X_i)=\det(X)^{r^{m-1}}
\end{equation}
for all $i$ and similarly for $Y_i, U_i$ and $V_i$ while $\det(W_i)=\det(W)^{r^{m-2}}$ for all $i$.

Now we first observe that $\det(X)=1$ since
\begin{equation*}
\sum_{i=1}^r (i-1)=\frac{r(r-1)}{2}
\end{equation*}
is divisible by $r$ due to $r$ being odd and hence $\det(X_i)=1$ for all $i$. Similarly, since the permutation corresponding to $Y$ is a cycle of odd length, we have $\det(Y)=1$ and thus $\det(Y_i)=1$. This means that the rescaling step can be omitted entirely for the matrices $X_i$ and $Y_i$.

In a similar spirit, note that
\begin{equation*}
\sum_{i=1}^r \frac{i(i-1)}{2}=\frac{(r-1)r(r+1)}{6}
\end{equation*}
is divisible by $r$ if $r>3$ and thus $\det(U)=1$ in this case, while the case $r=3$ yields $\det(U)=\omega$. Since we are assuming $d\geq 5$ and therefore $r=3$ implies $m\geq 2$, we get that $\det(U_i)=\det(U)^{r^{m-1}}=1$ for all $i$ in any case, so the matrices $U_i$ need not be rescaled either.

For the matrix $W$, we claim that $\det(W)=1$ and hence $\det(W_i)=1$ for all $i$ as well. Indeed, we show that the permutation $w$ defining $W$ is a product of $r-1$ cycles of length $r$, from which the result follows as $r$ is odd. To this end, take any $0\leq a\leq r^2-1$ and write it uniquely as $a=ur+v$ for some $0\leq u, v\leq r-1$. If $v>0$, then one computes
\begin{equation*}
\begin{split}
w(a)=(ur+v+(v-1)r)\dispmod r^2&=((u+v-1)r+v)\dispmod r^2 \\
&=((u+v-1)\dispmod r)\cdot r+v
\end{split}
\end{equation*}
and one easily verifies that this formula stays true if $v=0$. Thus, we see that applying $w$ leaves the value of $v$ invariant and maps $u$ to $(u+v-1)\inmod r$. We conclude that all $a$ with $v=1$ are fixed points of $w$; moreover, since $r$ is prime, any $a$ with $v\neq 1$ is part of a cycle 
\begin{equation*}
\begin{split}
v, \hspace{0.5cm} ((v-1)\dispmod r)\cdot r+v, \hspace{0.5cm} &(2(v-1)\dispmod r)\cdot r+v, \hspace{0.5cm}\dots
\end{split}
\end{equation*}
of length $r$ containing all $a$ having the same value of $v$. Hence, our claim is proved.

Thus, the only matrices for which the final rescaling step is necessary are the matrices $V_i$. Here, note that, to rescale one of the matrices $V_i\in\GL(d, q)$ for $1\leq i\leq m$ to have determinant $1$, we have to determine a $d$-th root $\xi$ of $\det(V_i)$ in $\F_q$, which, due to $d=r^m$ and $\det(V_i)=\det(V)^{r^{m-1}}$, may be done by determining an $r$-th root of $\det(V)$ in $\F_q$. Once this is done, multiplying each entry of $V_i$ with $1/\xi$ can be done in time $O(d^2\log^{1+\eps} q)$; as there are $m=O(\log d)$ matrices $V_i$, this has the required runtime.

To explain how to determine $\xi$, we begin with the case $r>3$ and examine $\det(V)$. Namely, since $V$ is a Vandermonde matrix, its determinant can be computed as
\begin{equation*}
\begin{split}
\det(V)&=\prod_{0\leq i<j\leq r-1} (\omega^j-\omega^i)=\prod_{i=0}^{r-2} \prod_{j=i+1}^{r-1} \omega^i(\omega^{j-i}-1)=\prod_{i=0}^{r-2} \prod_{k=1}^{r-1-i} \omega^i(\omega^k-1) \\
&=\prod_{i=0}^{r-2} \omega^{i(r-1-i)}\prod_{k=1}^{r-1-i} (\omega^k-1)=\prod_{i=0}^{r-2}\prod_{k=1}^{r-1-i} (\omega^k-1)\;,
\end{split}
\end{equation*}
where, in the last step, we have used that
\begin{equation*}
\sum_{i=0}^{r-2} i(r-1-i)=\frac{(r-2)(r-1)r}{6}
\end{equation*}
is divisible by $r$ due to $r>3$. Manipulating the above expression for $\det(V)$ further by switching the order of the products, we obtain
\begin{equation}
\label{eq:detV1}
\det(V)=\prod_{k=1}^{r-1}\prod_{i=0}^{r-1-k} (\omega^k-1)=\prod_{k=1}^{r-1} (\omega^k-1)^{r-k}\;.
\end{equation}
However, using that $\omega^r=1$, we can also write this as
\begin{equation*}
\det(V)=\prod_{k=1}^{r-1} (\omega^k-\omega^r)^{r-k}=\prod_{k=1}^{r-1}\omega^{k(r-k)}(1-\omega^{r-k})^{r-k}=\prod_{k=1}^{r-1} (1-\omega^{r-k})^{r-k}\;,
\end{equation*}
where the last step uses that
\begin{equation*}
\sum_{k=1}^{r-1} k(r-k)=\frac{(r-1)r(r+1)}{6}
\end{equation*}
is divisible by $r$ due to $r>3$. Reindexing the final product above then gives us
\begin{equation}
\label{eq:detV2}
\det(V)=\prod_{k=1}^{r-1} (1-\omega^k)^k=(-1)^{(r-1)/2}\prod_{k=1}^{r-1} (\omega^k-1)^k\;,
\end{equation}
where the factor $(-1)^{(r-1)/2}$ accounts for the additional sign we have incurred by replacing $1-\omega^k$ by $\omega^k-1$: indeed, if $r\inmod 4=1$, then there is an even number of odd integers $1\leq k\leq r-1$ and we do not need to account for an additional sign, while in the case $r\inmod 4=3$, there is an odd number of odd integers $1\leq k\leq r-1$ and we need to account for an additional minus sign. Finally, we multiply the expressions (\ref{eq:detV1}) and (\ref{eq:detV2}) to obtain
\begin{equation}
\begin{split}
\label{eq:detVsquare}
\det(V)^2&=(-1)^{(r-1)/2}\left(\prod_{k=1}^{r-1} (\omega^k-1)\right)^r \\
&=(-1)^{(r-1)/2}\left(\prod_{k=1}^{r-1} (1-\omega^k)\right)^r=(-1)^{(r-1)/2}r^r\;,
\end{split}
\end{equation}
where the last step uses the observation that the product occurring in the penultimate expression is just the $r$-th cyclotomic polynomial 
\begin{equation*}
\Phi_r(X)=\prod_{k=1}^{r-1} (X-\omega^k)=X^{r-1}+X^{r-2}+\dots+1
\end{equation*}
evaluated at $X=1$.

Thus, by (\ref{eq:detVsquare}), to obtain an $r$-th root of $\det(V)$ in $\F_q$, we may alternatively compute the square roots $\pm\xi'$ of $(-1)^{(r-1)/2}r$ in $\F_q$ (note that these square roots must exist in $\F_q$ since the $r$-th power of $(-1)^{(r-1)/2}r$ is a square in $\F_q$ by (\ref{eq:detVsquare})); using the big prime Berlekamp algorithm from \cite[Thm. 8.12]{GeddesCzaporLabahn}, this may be done in $O(\log^{2+\eps} q)$. Then at least one of $\pm\xi'$ must be an $r$-th root of $\det(V)$ and we may test this in time $O(r\log r\log^{1+\eps} q)$ by computing $\det(V)$ using the expression (\ref{eq:detV1}) and comparing the result with $\pm{\xi'}^r$.

It remains to deal with the case $r=3$. Here, one quickly computes $\det(V)=3\omega(\omega-1)$ and noting that $(1-\omega)^3=3\omega(\omega-1)$ then provides a third root of $\det(V)$ in $\F_q$, as required. Clearly, this can be computed in time $O(\log^{2+\eps} q)$.

\subsection{Maximal subgroups of type $(4\circ 2_-^{1+2m}).\Sp(2m, 2)$} 

Now let $E$ be a $2$-group of symplectic type; in particular, we have $r=2$. Recalling that we have assumed that $d=r^m>4$, we also have $m\geq 3$. By \cite[Lem. 9.3]{HoltRoneyDougal}, the definitions of the matrices $X_i, Y_i, V_i$ and $W_i$ stay the same as in the previous case. The only change we have to make is that $U\in\GL(2, q)$ is now the diagonal matrix with $U_{11}=1$ and $U_{22}=\zeta^{(q-1)/4}$, where $\zeta$ denotes a primitive element of $\F_q$, and the matrices $U_i$ for $1\leq i\leq m$ change accordingly.

Now one easily computes that $\det(X)=\det(Y)=\omega$ while $\det(U)=\zeta^{(q-1)/4}$ and $\det(W)=\pm 1$ since $W$ is a permutation matrix. Due to $m\geq 3$ and $r=2$, this means that $\det(X_i)=\det(Y_i)=\det(U_i)=\det(W_i)=1$ for all $i$ by (\ref{eq:detXi}) and the analogous relations for $Y_i, U_i$ and $W_i$. 

Thus, as in the previous case, only the matrices $V_i$ for $1\leq i\leq m$ need to be rescaled to determinant $1$ and for this, it again suffices to determine a $d$-th root of $\det(V_i)$ in $\F_q$ in the required time. To this end, one easily computes that $\det(V)=\omega-1=-2$. We now distinguish three cases:

\begin{enumerate}[label=(\roman*)]
\item $p\inmod 8=1$ or $p\inmod 8=3$. In this case, $-2$ is a square mod $p$ by quadratic reciprocity. Hence, we can find a square root of $-2$ in $\F_p$ and, due to $\det(V)=-2$, this will be a $d$-th root of $\det(V_i)=\det(V)^{2^{m-1}}$ for all $i$.

\item $p\inmod 8=7$. In this case, $2$ is a square mod $p$ by quadratic reciprocity. Hence, we can find a square root of $2$ in $\F_p$ and, due to $\det(V)^2=(-2)^2=2^2$, this will be a fourth root of $\det(V)^2$ and thus a $d$-th root of $\det(V_i)=(\det(V)^2)^{2^{m-2}}$ for all $i$.

\item $p\inmod 8=5$. In this case, $-1$ is a square mod $p$ by quadratic reciprocity, so we can find $j\in\F_p$ with $j^2=-1$. However, since $j^{(p-1)/2}=(-1)^{(p-1)/4}=-1$, by Euler's criterion, the element $j$ is not a square mod $p$ itself. As $2$ is not a square mod $p$ either by quadratic reciprocity, we see that we can find a square root of $2j$ in $\F_p$. Due to $\det(V)^4=(-2)^4=(2j)^4$, this will be an eighth root of $\det(V)^4$ and thus a $d$-th root of $\det(V_i)=(\det(V)^4)^{2^{m-3}}$ for all $i$.
\end{enumerate}

In all cases, the necessary square roots can be computed in time $O(\log^{2+\eps} q)$ by the big prime Berlekamp algorithm, as before, so we are done.

\subsection{Maximal subgroups of type $2_-^{1+2m}.O^-(2m, 2)$} 

Finally, let $E$ be an extraspecial $2$-group of minus type; in particular, we again have $r=2$ and $m\geq 3$. By \cite[Lem. 9.4]{HoltRoneyDougal}, the definitions of the matrices $X_i, Y_i, V_i$ and $W_i$ for $i>1$ stay the same as in the previous case. While the matrices $U_i$ for $i>1$ are to be discarded entirely, the matrices $X_1, Y_1, U_1, V_1$ and $W_1$ are now defined as
\begin{align*}
X_1&\coloneqq I_{2^{m-1}}\tensor\begin{pmatrix} a & b \\ b & -a\end{pmatrix}\;, \hspace{0.5cm} Y_1\coloneqq I_{2^{m-1}}\tensor\begin{pmatrix} 0 & -1 \\ 1 & 0\end{pmatrix}\;, \hspace{0.5cm} U_1\coloneqq I_{2^{m-1}}\tensor\begin{pmatrix} 1 & 1 \\ -1 & 1\end{pmatrix}\;, \\
V_1&\coloneqq I_{2^{m-1}}\tensor\begin{pmatrix} 1+a+b & 1-a+b \\ -1-a+b & 1-a-b\end{pmatrix}\;, \hspace{0.5cm} W_1\coloneqq I_{2^{m-2}}\tensor\begin{pmatrix} 1 & 0 & 1 & 0 \\ 0 & 1 & 0 & 1 \\ 0 & 1 & 0 & -1 \\ -1 & 0 & 1 & 0\end{pmatrix}\;,
\end{align*}
where $a, b\in\F_q$ are chosen such that $a^2+b^2=-1$.

As the other generators are the same as in the previous case, we only have to explain how to rescale $X_1, Y_1, U_1, V_1$ and $W_1$ to determinant $1$. Here, one easily calculates that $\det(X_1)=\det(Y_1)=1$, so no rescaling is necessary for $X_1$ and $Y_1$. Moreover, since $\det(V_1)=4^{2^{m-1}}=2^d$, the matrix $V_1$ merely has to be scaled by $1/2$. Finally, one calculates $\det(U_1)=\det(W_1)=2^{2^{m-1}}$; since $m\geq 3$, this is the same as $(-2)^{2^{m-1}}$ and thus computing a $d$-th root of $\det(U_1)=\det(W_1)$ can be done in exactly the same way as the computation of a $d$-th root of $\det(V_i)$ in the previous case. This finishes the proof.

\section*{Acknowledgements}
 
The author is grateful to Max Horn for introducing him to the topic of constructing maximal subgroups of classical groups and for help in preparing an earlier draft. The author also thanks Derek Holt for useful feedback on an earlier draft and Igor Shparlinski for his encouragement to publish these results.

\bibliographystyle{alpha}

\end{document}